\documentclass[12pt]{amsart}
\usepackage{amssymb}
\usepackage{amsmath}
\usepackage{amsthm}
\usepackage{amscd}

\renewcommand{\phi}{\varphi}

\newcommand{\prt}{\partial}
\newcommand{\E}{\mathop{\mathbf{E}_\mu}}
\newcommand{\Ent}{\mathop{\mathbf{Ent}_\mu}}

\newtheorem{theorem}{Theorem}
\newtheorem{prop}{Proposition}
\newtheorem{lem}{Lemma}

\theoremstyle{definition}

\theoremstyle{remark}
\newtheorem{rem}{Remark}
\newtheorem{ex}{Example}

\newcommand{\be}{\begin{equation}}
\newcommand{\ee}{\end{equation}}

\newcommand{\ba}{\begin{aligned}}
\newcommand{\ea}{\end{aligned}}

\def \eps{\varepsilon}

\def\Re{\mathbb{R}}

\def\1{1\!\!\hbox{{\rm I}}}
\begin{document}


\title{Lift zonoid order and functional inequalities}

\author{Alexei M. Kulik}
\address{Institute of Mathematics Nat. Acad. Sci. of Ukraine, Tereshchenkivska str. 3,  01601, Kyiv, Ukraine}
\email{kulik.alex.m@gmail.com}
\thanks{}

\author{Taras D. Tymoshkevych}
\address{National Taras Shevchenko University of Kyiv, Academician Glushkov pr. 4-e, 03127, Kyiv, Ukraine}
\email{tymoshkevych@gmail.com}
\thanks{T.Tymoshkevych was partially supported by the Leonard Euler program, DAAD project No. 55518603.
}

\dedicatory{} 

\begin{abstract}
We introduce the notion of  a weighted lift zonoid and show that,
for properly chosen weights $v$, the ordering condition on a measure
$\mu$, formulated in terms of the weighted lift zonoids of this
measure, leads to certain functional inequalities for this measure,
such as non-linear extensions of Bobkov's {shift inequality} and
weighted {inverse log-Sobolev inequality}. The choice of the weight
$K$, involved in our version of the inverse log-Sobolev inequality,
differs substantially from those available in the literature, and
requires the weight $v$, involved into the definition of the
weighted lift zonoid, to equal the {divergence} of the weight $K$
w.r.t. initial measure $\mu$. We observe that such a choice may be
useful  for proving direct log-Sobolev inequality, as well, either
in its weighted or classical forms.
\end{abstract}

\maketitle

\section {Introduction}

The notions of \emph{zonoid} and \emph{lift zonoid}, introduced in
\cite{KM_AnnSta}, have a diverse field of  applications. Because the
lift zonoid determines the underlying measure uniquely, this
concept can be used in multivariate statistics for measuring the
variability of  laws of random vectors, and for ordering these laws,
see \cite{KM_Be}. The concept of \emph{zonoid equivalence} appears
to be both naturally motivated by financial applications, and useful
for proving extensions of the ergodic theorem for zonoid stationary
and zonoid swap-invariant random sequences, see \cite{MS11},
\cite{MSS12}.  Lift zonoids lead   naturally to definitions of
associated $\alpha$-trimming and data depth, see \cite{KM_AnnSta}
and \cite{Cascos2010}, and to barycentric representation of the
points of a space with a given measure, see \cite{KM_AnnSta} and
\cite{KulTym}.

In this paper, we explore a new field, where the notion of lift
zonoid can be applied naturally. As a straightforward extension of
the definition of lift zonoid, we introduce a \emph{weighted lift
zonoid} $\hat Z^v(\mu)$ with a vector-valued weight function $v$. We
show that, for properly chosen weights $v$, the ordering condition
on a measure $\mu$, formulated in terms of the weighted lift zonoid
of this measure, leads to certain functional inequalities for this
measure, such as non-linear extensions of Bobkov's \emph{shift
inequality} \cite{Bo} and weighted \emph{inverse log-Sobolev
inequality}.  Weighted versions of the classical functional
inequalities (Poincar\'e, log-Sobolev, etc) have been studied
recently in various contexts. The choice of the weight $K$, involved
in our version of inverse LSI, is specific and differs substantially
from those available in the literature. This choice is strongly
motivated by (an extension of) the functional form of Bobkov's shift
inequality, and requires the weight $v$, involved into the
definition of the weighted lift zonoid, to equal the
\emph{divergence} of the weight $K$ w.r.t. initial measure $\mu$. We
observe that such a choice may be useful  for proving (weighted)
direct log-Sobolev inequality, as well. In the case of a bounded
weight, this may lead to new sufficient conditions for the
log-Sobolev inequality. We illustrate the range of applications of
these conditions in two examples in Section \ref{s4} below.

\section{Weighted lift zonoids, non-linear shift inequalities, and  weighted inverse log-Sobolev inequalities}\label{s:1}
Let $\mu$ be a probability measure on the Borel $\sigma$-algebra in
$\Re^d$, and $v:\Re^d\to \Re^d$ be a measurable function such that
$$
\int_{\Re^d}\|v(x)\|\mu(dx)<\infty;
$$
here and below we denote by $\|\cdot\|$ the Euclidean norm in
$\Re^d$. We define the \emph{weighted zonoid}  $Z^v(\mu)$ with the
weight $v$ as the set of all the points in $\Re^d$ of the form
\be\label{int_mu} \int_{\Re^d}g(x)v(x)\, \mu(dx) \ee with arbitrary
Borel measurable $g:\Re^d\to [0,1]$.  The \emph{weighted lift
zonoid} $\hat Z^v(\mu)$ is defined as the weighted zonoid of the
measure $\delta_1\times\mu$ in $\Re^{d+1}$. Equivalently, the
weighted zonoid $Z^v(\mu)$ and the weighted lift zonoid $\hat
Z^v(\mu)$ are the sets of the points of the form \be\label{expect}
Eg(X)v(X)\in \Re^d\quad \hbox{and}\quad \Big(Eg(X),
Eg(X)v(X)\Big)\in \Re^{d+1} \ee respectively, where $X$ is a random
vector with the distribution $\mu$.  This definition is a
straightforward generalization of the definitions of the zonoid and
the lift zonoid (see \cite{KM_Be}, Definition 2.1), where the
function $v$  has the form $v(x)=x$.

The  lift zonoid  $\hat Z(\mu)$ is a  convex compact set in
$\Re^{d+1}$, symmetric w.r.t. the point  $((1/2),(1/2)EX)$, which
identifies the underlying measure $\mu$ uniquely; see \cite{KM_Be}.
On the other hand, it can be seen easily that the definition of the
weighted lift zonoid $\hat Z^v(\mu)$ would not change if  one
restricts the class of Borel measurable functions $g$ within it to
the class of the functions of the form
$$
g(x)=G(v(x)), \quad  G:\Re^d\to [0,1]\hbox{ is Borel measurable}.
$$
This observation leads immediately to the identity $\hat
Z^v(\mu)=\hat Z(\mu\circ v^{-1})$; that is, the weighted lift zonoid
$\hat Z^v(\mu)$ equals the (usual) lift zonoid of the image of the
measure $\mu$ under the mapping $v$. As a corollary, we get that the
weighted lift zonoid $\hat Z^v(\mu)$ is a  convex compact set in
$\Re^{d+1}$, symmetric w.r.t. the point  $((1/2),(1/2)Ev(X))$, and
identifies the image measure $\mu\circ v^{-1}$ uniquely.

The following theorem motivates the above definition of the weighted lift
zonoid. To formulate it, we need to introduce some
notation. Denote by $\gamma_c$ the centered Gaussian measure in
$\Re^d$ with the covariance matrix $c^2I_{\Re^d}$. Let
$$\phi(x)={1\over \sqrt{2\pi}} e^{-x^2/2}, \quad \Phi(x)=\int_{-\infty}^x\phi(y)\, dy, \quad x\in \Re$$
be the standard Gaussian distribution density function and the
standard Gaussian cumulative distribution function, respectively,
and let \be\label{iso} I(p)=\phi(\Phi^{-1}(p)), \quad p\in (0,1),
\quad I(0)=I(1)=0\ee be the \emph{Gaussian isoperimetric function}.

For any measurable $f$ on $\Re^d$, we write $\E f$ for the integral
of $f$ w.r.t. $\mu$; function $f$ may be vector-valued, then the
integral is understood in the component-wise sense. For a function
$f$ taking values in $\Re^+$, its $\mu$-entropy is defined by
$$
\Ent f=\E (f\log f)-(\E f)\log(\E f),
$$
with the convention $0\log 0=0$.

In what follows, we  assume that the measure $\mu$ has \emph{the
logarithmic gradient} $v_\mu$; that is, a function $v_\mu:\Re^d\to
\Re$, integrable w.r.t. $\mu$ and such that for every smooth
$f:\Re^d\to \Re$ with a compact support \be\label{IBP} \E\nabla
f=-\E (v_\mu f). \ee This assumption is equivalent to the following,
see Proposition 3.4.3 in \cite{Bogachev}: there exists the density
$p_\mu$ of the measure $\mu$ w.r.t. the Lebesgue measure, which
belongs to the Sobolev class $W_{1,1}(\Re^d)$; in this case
$$
[v_\mu]_i={\prt_{x_i}p_\mu\over p_\mu}, \quad i=1, \dots, d.
$$

\begin{theorem}\label{t1} I. The following three statements are equivalent.
\begin{itemize}
  \item[\textbf{A.}] $\widehat Z^{v_{\mu}}(\mu) \subset \widehat
Z(\gamma_c)$.
  \item[\textbf{B.}] For any smooth function $f: \Re^{d}
\rightarrow [0,1]$ with a compact support, one has
\be\label{shift_inf} \|\E \nabla f\| \leq c I (\E f).\ee
  \item[\textbf{C.}] For any $h\in \Re^{d}, A \in
\mathcal{B}(\Re^{d})$ \be\label{shift}
\Phi\Big(\Phi^{-1}(\mu(A))-c\|h\|\Big) \leq \mu(A+h) \leq
\Phi\Big(\Phi^{-1}(\mu(A))+c\|h\|\Big). \ee

\end{itemize}

II. Under the conditions \textbf{A - C} above, the following inverse
log-Sobolev inequality holds true: for any smooth function $f:
\Re^{d} \rightarrow [0,\infty)$ with a compact support,
\be\label{ILSE}\|\E \nabla f\|^2 \leq 2c\Ent f \E f.\ee

\end{theorem}

\begin{rem}\label{r1} By the definition (see Definition 5.1 in \cite{KM_Be}), two measures  $\mu_1, \mu_2$ are related by the \emph{lift zonoid order} (notation: $\mu_1\preccurlyeq_{LZ}\mu_2$), if
$$
\hat Z(\mu_1)\subset \hat Z(\mu_2).
$$
 Recall that  $\widehat Z^{v_\mu}(\mu)$ equals the lift zonoid of $\nu_\mu:=\mu\circ v_\mu^{-1}$; that is, of the distribution of the logarithmic gradient of the measure  $\mu$. Hence statement \textbf{A} can be equivalently formulated as follows: the distribution $\nu_\mu$ of the logarithmic gradient of the measure  $\mu$ is dominated in the sense of the lift zonoid order by the canonical Gaussian measure in $\Re^d$.
\end{rem}

Theorem \ref{t1} is not a genuinely new one. The equivalence of the
relations \textbf{B} and \textbf{C} is used by S. Bobkov in
\cite{Bo} as a key ingredient in the proof of the \emph{shift
inequality} (\ref{shift}) (in \cite{Bo}, the measure $\mu$ is
supposed to be a product-measure, but the proof of the equivalence
of  (\ref{shift_inf}) and (\ref{shift}) in fact does not rely on
this assumption). The outline of the proof of (\ref{ILSE}) under
(\ref{shift_inf}) and (\ref{shift}) is given in \cite{BCF}. What we
would like to emphasize is that condition \textbf{B}, usually called
the \emph{functional version of the shift inequality}, is equivalent
to the relation \textbf{A}, which according to Remark \ref{r1} can
be  written as the lift zonoid order relation \be\label{LZ} \nu_\mu
\preccurlyeq_{LZ}\gamma_c. \ee It is instructive to compare
(\ref{LZ}) with the following necessary and sufficient condition for
the functional version of the shift inequality to hold, given in
\cite{Bo} in the case where the measure $\mu$ is a product-measure
with equal marginals $\mu_1$. This condition states that there
exists $c>0$ such that (\ref{shift_inf}) holds true,  if and only if
there exists $\eps>0$ such that \be\label{prod_moment} \int_\Re
e^{\eps x^2}\nu_{\mu_1}(dx)\leq 2; \ee in addition, the optimal
constant $c$ in (\ref{shift_inf}) and $\eps$ in (\ref{prod_moment}) are connected
by the relation \be\label{rel} {1\over \sqrt{6\eps}}\leq c\leq
{4\over \sqrt{\eps}}. \ee For the product measure
$\mu(dx)=\prod_{i=1}^d\mu_1(dx_i)$, respective distribution  of the
logarithmic gradient is again a product measure
$$
\nu_\mu(dx)=\prod_{i=1}^d\nu_{\mu_1}(dx_i),
$$
and in this case, due to Corollary 5.3 in \cite{KM_Be},  (\ref{LZ})
is equivalent to \be\label{LZ_1} \nu_{\mu_1}
\preccurlyeq_{LZ}\gamma_c^1, \ee where $\gamma_c^1$ is the
$\mathcal{N}(0, c^2)$-Gaussian measure on $\Re$. Both
(\ref{prod_moment}) and (\ref{LZ_1}) are conditions on the tails of
the distribution  of the  logarithmic gradient of $\mu_1$, but
(\ref{LZ_1}) is more precise because it involves the same $c$ with
(\ref{shift_inf}).

The main result of this section, Theorem \ref{t2} below, is a generalization of Theorem \ref{t1} and is motivated by an observation that in Theorem \ref{t1} the equivalence of the relations \textbf{A} and \textbf{B} follows in a very straightforward way from the integration-by-parts formula (\ref{IBP}); see the proof of Theorem \ref{t2} below. With this observation in mind, we introduce a wide class of weights which admit an analogue of the integration-by-parts formula (\ref{IBP}). To do
that, we recall that the $\mu$-divergence of a function $g:\Re^d\to
\Re^d$, if exists, is defined as the function $\delta_\mu(g)\in
L_1(\Re^d, \mu)$ such that for every smooth $f:\Re^d\to \Re$ with a
compact support
$$
\E (\nabla f, g)_{\Re^d}=\E f\delta_\mu(g).
$$
The $\mu$-divergence is well defined, for instance, for any $g\in
C^1$ bounded together with its partial derivatives; in this case,
$$
\delta_\mu(g)=-\sum_{i=1}^d[v_\mu]_ig_i-\sum_{i=1}^d\prt_{x_i}g_i.
$$
This follows directly from (\ref{IBP}); see \cite{Bogachev}, Chapter
6 for more information on this subject. Let  function $v:\Re^d\to
\Re$ be such that,  for some function  $K$ taking values in $d\times
d$-matrices, \be\label{v} v_i=\delta_\mu(K_i), \quad i=1, \dots, d,
\ee where $K_i$ denotes the $i$-th row of the matrix $K$. Then for
every smooth $f$ with  a compact support \be\label{IBP_K} \E \Big(K
\nabla f\Big)=\E fv; \ee here and below we treat elements of $\Re^d$
as vectors-columns. Formula (\ref{IBP_K}) is a straightforward
extension of the integration-by-parts formula (\ref{IBP}), where
the gradient $\nabla$  is replaced by the ``weighted gradient''
$K\nabla$ with the matrix-valued weight $K$, and the logarithmic
gradient $v_\mu$  is replaced by the $\mu$-divergence of $K$.
Furthermore, if $K$ satisfies some extra regularity condition, e.g.
\be\label{K_reg} K:\Re^d\to \Re^{d\times d}\quad \hbox{is
Lipschitz}, \ee then for every $h\in \Re^d$ there exists a flow of
solutions $\{\Psi^{K,h}_t(x), t\in \Re, x\in \Re^d\}$ of the Cauchy
problem \be\label{ODE} d\Psi_t(x)=(K^*h)\Big(\Psi_t(x)\Big)\, dt,
\quad \Psi_0(x)=x. \ee

\begin{theorem}\label{t2} I. Let $v=(v_i)_{i=1}^d$ satisfy (\ref{v}). Then the following two statements are equivalent.
\begin{itemize}
  \item[\textbf{A1.}] $\widehat Z^{v}(\mu) \subset \widehat
Z(\gamma_c)$.
  \item[\textbf{B1.}] For any smooth function $f: \Re^{d}
\rightarrow [0,1]$ with a compact support,  one has
\be\label{shift_inf_K} \|\E K\nabla f\| \leq c I (\E f).\ee
\end{itemize}

If, in addition,  the matrix-valued function $K$ satisfies
(\ref{K_reg}), then \textbf{A1} and \textbf{B1} are equivalent to
the following.
\begin{itemize}
  \item[\textbf{C1.}] For any $h\in \Re^{d}, A \in
\mathcal{B}(\Re^{d})$ \be\label{shift_K}
\Phi\Big(\Phi^{-1}(\mu(A))-c\|h\|\Big) \leq
\mu\Big(\left[\Psi^{K,h}_1\right]^{-1}(A)\Big)  \leq
\Phi\Big(\Phi^{-1}(\mu(A))+c\|h\|\Big).
 \ee
\end{itemize}

II. Under the condition \textbf{A1}, equivalently \textbf{B1}, the
following weighted inverse log-Sobolev inequality holds true: for
any smooth function $f: \Re^{d} \rightarrow [0,\infty)$ with a
compact support, \be\label{ILSE_K}\|\E K\nabla f\|^2 \leq 2c^2\Ent f
\E f.\ee

\end{theorem}

Note that condition \textbf{A1} is just the lift zonoid order relation for the
image measure of $\mu$ under $v$: \be\label{LZ_v} \mu\circ v^{-1}
\preccurlyeq_{LZ}\gamma_c. \ee

Before giving the proof of Theorem \ref{t2}, let us summarize: a lift zonoid order condition (\ref{LZ}) is a
\emph{criterion} for the shift inequality, written  either in its
explicit form (\ref{shift}), or in its functional form
(\ref{shift_inf}). This equivalence is rather flexible in the
following sense: if the logarithmic gradient  $v_\mu$ in (\ref{LZ})
is replaced by another weight $v$ of the form \be\label{div}
v=\delta_\mu(K) \ee (see (\ref{v})), then respective lift zonoid
order condition (\ref{LZ_v}) is still equivalent to the weighted
version (\ref{shift_inf_K}) of the functional form of a
(generalized) shift inequality. The explicit form of the
(generalised)  shift inequality in that case is available as well,
and concerns, instead of linear shifts,  the transformations of the
initial measure $\mu$  by the flows of solutions to (\ref{ODE}).

{\it Proof of Theorem \ref{t2}: statement I.} The lift zonoid $\hat Z(\gamma)$ of a standard Gaussian measure
$\gamma$ in $\Re^d$ can be identified in the following way: for a
given $\alpha\in (0,1)$, the section of $\hat Z(\gamma)$ by the
hyper-plane $\{\alpha\}\times \Re^d$ has the projection on the last
$d$ coordinates equal to the ball centered at $0$ and having the
radius $I(\alpha)$; see \cite{KM_AnnSta}, Section 6.3 or
\cite{KulTym}, Proposition 3.4. It is easy to see from the
definition of the lift zonoid that
$$
\hat Z(\gamma_c)=c\hat Z(\gamma).
$$
Hence condition \textbf{A1} can be equivalently written as follows:
for every Borel measurable $g:\Re^d\to [0,1]$ such that $ \E
g=\alpha, $
$$
\left\|\E \Big(g v\Big)\right\|\leq cI(\alpha)=cI(\E g).
$$
By the standard approximation argument, the above condition is
equivalent to a similar one with  Borel measurable $g$'s replaced by
smooth and compactly supported $f$'s. Because for such $f$ by
(\ref{IBP_K})
$$
\left\|\E \Big(f v\Big)\right\|=\left\|\E \Big(K\nabla
f\Big)\right\|,
$$
conditions \textbf{A1} and \textbf{B1} are equivalent.

The proof of the equivalence of \textbf{B1} and \textbf{C1} follows
the same lines with the S.Bobkov's proof from \cite{Bo} for the case
of product measures and linear shifts; to make the exposition
self-sufficient here we expose the key steps of this proof.

Denote $R_{r}(p)=\Phi\Big(\Phi^{-1}(p)+r\Big), r\geq 0, p\in (0,1)$.
Then the following properties hold true:
\begin{itemize}
  \item for every $r\geq 0$ the function $R_r$ is concave;
  \item the family $\{R_r, r\geq 0\}$ is a semigroup w.r.t. the composition of the functions, i.e.
  $$
  R_{r_1}\circ R_{r_2}=R_{r_1+r_2};
  $$
  \item the function $R_0$ is an identity, and the ``generator'' of the semigroup $\{R_r, r\geq 0\}$ equals the Gaussian isoperimetric function $I$ in the sense that
  $$
  {R_r(p)-p\over r}\to I(p), \quad r\to 0+.
  $$
\end{itemize}
Similarly, the family of functions
$S_{r}(p)=\Phi\Big(\Phi^{-1}(p)-r\Big), r\geq 0, p\in (0,1)$ has the
following properties:
\begin{itemize}
  \item for every $r\geq 0$ the function $S_r$ is convex;
  \item the family $\{S_r, r\geq 0\}$ is a semigroup w.r.t. the composition of the functions;
  \item the function $S_0$ is an identity, and the ``generator'' of the semigroup $\{S_r, r\geq 0\}$ equals  $(-I)$.
\end{itemize}

Observe that \textbf{C1} is equivalent to the following.

\begin{itemize}
  \item[\textbf{C2.}] For any $h\in \Re^{d}$ and Borel measurable $f:\Re^d\to [0,1]$
\be\label{shift_K_f} S_{c\|h\|}(\E f) \leq  \E
\Big(f\circ\Psi^{K,h}_1\Big)\leq   R_{c\|h\|}(\E f).
 \ee
\end{itemize}

Indeed, taking $f=\1_A$ we get \textbf{C1} from \textbf{C2}.
Inversely, under \textbf{C1} by the concavity of $R_r$ and Jensen's
inequality we have
$$\ba
\E \Big(f\circ\Psi^{K,h}_1\Big)&=
\int_{0}^\infty\mu(\{x: f(\Psi^{K,h}_1(x))\geq t\})\, dt\leq \int_{0}^\infty R_{c\|h\|}\Big(\mu(\{x: f(x)\geq t\})\Big)\, dt\\
&\leq R_{c\|h\|}\left(\int_{0}^\infty \mu(\{x: f(x)\geq t\})\,
dt\right)= R_{c\|h\|}(\E f). \ea
$$
The proof of the left hand side inequality in (\ref{shift_K_f}) is
similar and omitted. Hence \textbf{C1} and  \textbf{C2} are
equivalent.

To get \textbf{B1} from \textbf{C2}, take $th$ instead of $h$ and
differentiate the right hand side inequality in  (\ref{shift_K_f})
w.r.t. $t$ at the point $t=0$. In  more details, denote $f_t(x)=
f(\Psi^{K,th}_1(x))$, then
$$
f_t(x)=f(\Psi^{K,h}_t(x)),
$$
and therefore there exits a continuous derifative
$$
\prt_tf_t(x)=\Big((\nabla f) (\Psi^{K,h}_t(x)), (K^*h)
(\Psi^{K,h}_t(x))\Big)_{\Re^d}.
$$
Because $f$ is smooth and compactly supported and $K$ satisfies
(\ref{K_reg}), this derivative is bounded as a function of $(t,x)\in
[0, T]\times\Re^d$ for every fixed $T$. Therefore by the dominated
convergence theorem
$$
{1\over t}\Big(\E f_t -\E f\Big)\to \E \Big(\nabla f,
K^*h\Big)_{\Re^d}=\Big(\E K\nabla f,h\Big)_{\Re^d}, \quad t\to 0+.
$$
Because
$$
{1\over t}\Big(R_{ct\|h\|}(\E f)-\E f\Big)=c\|h\|I(\E f),
$$
we get from (\ref{shift_K_f})
$$
\Big(\E K\nabla f,h\Big)_{\Re^d}\leq c\|h\|I(\E f), \quad h\in
\Re^d.
$$
Taking $\sup$ over all $h$ with $\|h\|=1$, we get
(\ref{shift_inf_K}).

To get  \textbf{C2} from   \textbf{B1}, consider first the case
where $f$ is smooth and compactly supported and such that $0<\E
f<1$. By (\ref{shift_inf_K}), for a given $h\in \Re^d$ we have that
$$
{d\over dt}\Big|_{t=0}\E f_t=\Big(\E K\nabla f,h\Big)_{\Re^d}\leq
c\|h\|I(\E f).
$$
Recall that
$$
{d\over dt}\Big|_{t=0}R_{c\|h\|t}(\E f)=c\|h\|I(\E f).
$$
Therefore for every $\varrho>1$ there exists $\delta=\delta(f)>0$
such that for every $t\in (0,\delta)$: \be\label{31}\E f_t \leq
R_{\varrho c\|h\|t}(\E f). \ee Note that if $t_1\in (0, \delta(f))$
and $t_2\in (0, \delta(f_{t_1}))$, then \be\label{32} \E
f_{t_1+t_2}=\E \Big(f_{t_1}\circ\Psi^{K,h}_{t_2}\Big)\leq R_{\varrho
c\|h\|t_1}(\E f_{t_1})\leq R_{\varrho c\|h\|(t_1+t_2)}(\E f); \ee
here we have used the flow property of $\{\Psi^{K,h}_{t}, t\in
\Re\}$, the semigroup property of $\{R_r, r\geq 0\}$, and
monotonicity of $R_r$. Because the derivative $\prt_t f_t$ is
uniformly continuous w.r.t. $(t,x)\in [0, T]\times\Re^d$ for every
fixed $T$, it can be shown that
$$
\delta_T=\inf_{t\in [0,T]}\delta(f_t)>0.
$$
Then, applying (\ref{32}) at most $T/\delta_T$ times,  we get that
(\ref{31}) holds true for every $t\in [0, T]$. Consequently,
(\ref{31}) holds true for every $t\in \Re^+$ and $\varrho$ therein
can be replaced by 1. This gives the right hand side inequality in
(\ref{shift_K_f}) for  smooth and compactly supported $f$ such that
$0<\E f<1$. By an approximation argument, this can be extended to
any measurable $f:\Re^d\to [0,1]$. The proof of the left hand side
inequality in (\ref{shift_K_f}) is completely analogous and omitted.\qed

{\it Proof of Theorem \ref{t2}: statement II.} The following lemma is a straightforward extension of a part of
Proposition 2 in \cite{BCF} (the one which states the equivalence of
$P_1(c)$ and $P_2(c)$ in the notation of \cite{BCF}).

\begin{lem}\label{l1} Statement \textbf{B1} is equivalent to the following.

\begin{description}
   \item[B2.] For any smooth function $f: \Re^{d}
\rightarrow [0,1]$ with a compact support, one has
\be\label{shift_inf_k_iso} \sqrt{(\E I(f))^{2}+\frac{1}{c^2}\|\E
K\nabla f\|^{2}} \leq I (\E f). \ee
\end{description}
\end{lem}
{\it The proof } is completely analogous to the one from \cite{BCF},
therefore we just sketch it. The implication \textbf{B2}
$\Rightarrow$ \textbf{B1} is trivial. To get the inverse
implication, recall first that the standard Gaussian measure
$\gamma^d$ on $\Re^d$ satisfies \textbf{B2} with $c=1$ and identity
matrix $K$; see \cite{BCF}, Section 2. Consider a smooth function
$f:R^{d} \rightarrow [0,1]$ with a compact support, and let
$F(r)=\mu(\{x: f(x) \leq r\})$ be its distribution function w.r.t.
$\mu$. Assume that $F$ is absolutely continuous w.r.t. Lebesque
measure on $\Re$, and take $r \in \Re$, $\varepsilon
>0$. Define $\psi_{\varepsilon} (x) = \mathbb{I}_{[0,r]}(x)+(1-\frac{x-r}{\varepsilon})\mathbb{I}_{[r,r+\varepsilon]}(x).$
Applying \textbf{B1} to the function $g=\psi_{\varepsilon}(f)$ and
tending $\varepsilon \rightarrow 0$, we get \be\label{41}
F^{'}(r)\|\theta(r)\| \leq c I(F(r))\quad\hbox{for}\quad \mu\circ
f^{-1}\hbox{-a.a. }r\in \Re,  \ee where  $\theta(r)=E_{\mu}(K\nabla
f| f=r)$. Denote $k=F^{-1}\circ \Phi$, then $k$ transforms the
standard Gaussian measure $\gamma^1$ on $\Re$ to $\mu\circ f^{-1}$.
Taking the derivative in the identity $F(k)=\Phi,$  we get
$k^{'}F^{'}(k) =\phi$. Then from (\ref{41}) with $r=k(x)$ we get
inequality \be\label{42}\frac{1}{c}\|\theta(k(x))\| \leq k^{'}(x)\ee
valid $\gamma^1$-a.s. We have already mentioned that a standard
Gaussian measure satisfies \textbf{B2} with $c=1$ and identity $K$;
for the case $d=1$ this can be written as
$$ \sqrt{\left(\int_{\Re} I(g)d\gamma^{1}\right)^{2}+\left(\int_{\Re} g^{'}d \gamma^{1}\right)^{2}} \leq I\left(\int_{\Re} g d \gamma^{1}\right).$$
Applying this inequality to $g=k$ and using (\ref{42}) we get
\be\label{43} \sqrt{\left(\int_0^1 I(r)d
F(r)\right)^{2}+\left(\int_{0}^{1} \frac{1}{c}\|\theta(r)\|d
F(r)\right)^{2}} \leq I\left(\int_{0}^{1} r d F(r)\right);  \ee here
we have took into account that the image of $\gamma^1$ under $k$ is
$\mu\circ f^{-1}$, and $\mu\circ f^{-1}$ is supported in $[0,1]$.
Using the inequality $$ \int_0^1 \|\theta(r)\| d F(r) \geq
\left\|\int_0^1 \theta(r) d F(r)\right\|=\|\E K\nabla f \|,$$ we
complete the proof of the required statement. The additional
assumption of $\mu\circ f^{-1}$ to be absolutely continuous can be
removed by an approximation argument. \qed

According to Lemma \ref{l1}, to prove statement II of Theorem
\ref{t2} it is enough to show that \textbf{B2} implies
(\ref{ILSE_K}) for any non-negative smooth compactly supported $f$.
Take $\eps$ small, then  $\eps f$ takes values in $[0,1]$ and one
can apply \textbf{B2}. After trivial transformations, we get
$$\frac{1}{c^2}\|\E K\nabla f \|^{2} \leq \frac{I^{2}(\eps\E f)-(\E I(\varepsilon f))^{2}}{\varepsilon^{2}}.$$
Hence the required statement would follow from the relation
\be\label{44} \lim_{\varepsilon \rightarrow 0+} \frac{I^{2}(\eps\E
f)-(\E I(\varepsilon f))^{2}}{\varepsilon^{2}}=2\Ent f \E f. \ee
This relation can be proved straightforwardly using the following
asymptotic expansion: \be\label{45}I(\varepsilon)=\varepsilon\sqrt{2
\log\frac{1}{\varepsilon}}-\frac{\varepsilon\log (2\log
\frac{1}{\varepsilon})}{2\sqrt{2 \log
\frac{1}{\varepsilon}}}+\frac{\varepsilon}{\sqrt{2 \log
\frac{1}{\varepsilon}}}+\frac{\varepsilon\kappa(\varepsilon)}{\sqrt{2\log
\frac{1}{\varepsilon}}},\ee where $\kappa(\varepsilon) \rightarrow
0, \varepsilon \rightarrow 0 +$; the detailed exposition is
straightforward but cumbersome and therefore is omitted. The
asymptotic expansion (\ref{45}) follows from the standard expansion
$$ \Phi (t)= -\frac{1}{t}\phi(t)+\frac{1}{t^{3}}\phi(t)+O(t^{-5}\phi(t)),\quad  t \rightarrow -\infty,$$
which holds true e.g. by the integration-by-parts formula.\qed
\begin{rem} The above proof of statement II follows, in main lines, the one sketched in \cite{BCF} (the proof of the implication $P_3(c)\Rightarrow P_6(c\sqrt{2})$ in Proposition 2), where the authors referred to Beckner's lectures at the Institut Henri Poincar\'e. However, instead of using  the equivalence $I(\varepsilon)\sim\varepsilon\sqrt{2
\log\frac{1}{\varepsilon}},  \eps\to 0$, which apparently is not sufficient to provide (\ref{44}), we use stronger asymptotic expansion (\ref{45}).
\end{rem}

At the end of this section, let us mention that a more explicit condition, sufficient for the lift zonoid relation  (\ref{LZ_v}) tohold true,  can be given in a way similar to (\ref{prod_moment}).

\begin{prop}\label{p1} There exists $c>0$ such that (\ref{LZ_v}) holds true,  if and only if, there exists $\eps>0$ such that
\be\label{prod_moment_v} \E e^{\eps (v,h)^2_{\Re^d}}\leq 2, \quad
\|h\|\leq 1. \ee The optimal constant $c$ in (\ref{LZ_v}) and $\eps$ in
(\ref{prod_moment_v}) are connected by the relation (\ref{rel}).
\end{prop}

Because the lift zonoid order relation is equivalent to the same relation for all one-dimensional projections (see Section 5 in  \cite{KM_Be}), statement of  Proposition \ref{p1} follow immediately from the one-dimensional statement given below.

\begin{lem}\label{lB} For a measure $\nu$ on $\Re$  there exists $c>0$ such that
$$
\nu\preccurlyeq_{LZ}\gamma_c
$$
with  $\gamma_c\sim\mathcal{N}(0, c)$ if, and only if, there exists $\eps>0$ such that
$$
\int_{\Re}e^{\eps x^2}\, dx\leq 2;
$$
in that case, the optimal constants $c, \eps$ are  connected by the relation (\ref{rel}).
\end{lem}

The proof of Lemma \ref{lB} is contained, in fact, in the proof of Lemma 4.1 in \cite{Bo}, hence we omit it here.

\section{Weighted log-Sobolev inequalities in $\Re$}\label{s:1}

Theorem \ref{t2} above gives a sufficient condition for a
weighted inverse log-Sobolev inequality, based on a pair of
functions $v, K$ related by (\ref{div}).  The main result of this section, Theorem \ref{t3} below,
shows
that the use of  the same pair may lead to  sufficient conditions
for the (direct) log-Sobolev inequality, either in a weighted or in
a classical form. What is surprising is that, even in the  simplest
one-dimensional case, Theorem \ref{t3} leads to new sufficient
conditions for the log-Sobolev inequality, when compared with those
available in a literature; see  below Proposition \ref{p2}, Proposition \ref{p3}, and two
examples in Section \ref{s4}. We believe that the reason for
that is a proper choice of the \emph{pair} of the weight functions
$v,K$, involved in (\ref{BE_K}) and connected by (\ref{div}).

\begin{theorem}\label{t3} Let $d=1$ and functions $v, K$ be related by (\ref{div}). Assume that for some $\alpha>0$
\be\label{BE_K}
 Kv^{'} \geq \alpha.
\ee Assume in addition that the functions $K$ and \be\label{ab}
a:=2KK'+K^2 v_\mu \ee belong to $C^\infty$, have at most linear
growth at $\infty$,  and all their derivatives have at most
polynomial growth at $\infty$.

Then for every smooth $f$ with a compact support \be\label{LSI_K}
\Ent f^{2} \leq \frac{2}{\alpha}\E\Big(K f^{'}\Big)^2.\ee

As a corollary, if $K$ is bounded then $\mu$ satisfies the
(classical) log-Sobolev inequality: for every absolutely continuous
$f$ such that both $f$ and $f'$ are square integrable w.r.t. $\mu$,
\be\label{LSI} \Ent f^{2} \leq
\frac{2}{\alpha}\Big(\sup_xK^2(x)\Big)\E\Big(f^{'}\Big)^2. \ee

\end{theorem}

\begin{rem}
The proof of Theorem \ref{t3} is based on the classic Bakry-Emery
criterion; see below. We strongly
believe that similar technique  is applicable in the multidimensional case
as well, but because of possible non-commutativity of matrix-valued
weights which appear therein, now we can not give a multidimensional
version of Theorem \ref{t3}; this is a subject for a further
research.
\end{rem}

\begin{rem} The additional assumptions on the functions $K, a$ to be smooth and to satisfy certain growth bounds, in particular cases, can be removed by an approximation procedure; see e.g. Propositions \ref{p2} and \ref{p3} below.
\end{rem}

{\it Proof of Theorem \ref{t3}.} Consider a Markov process $X$ defined as the strong solution to the
SDE
$$
dX_t=a(X_t)\, dt+\sqrt{2}K(X_t)\, dW_t;
$$
see (\ref{ab}) for the formula for the coefficient $a$. Then on the
Schwartz space $\mathcal{S}(\Re)$ of  $C^\infty$ functions s.t. all
their derivatives decay at $\infty$ faster than any polynomial, the
generator $L$ of the process $X$ has the form
$$
Lf=af'+bf''=v_\mu f'+(bf')', \quad b:=K^2.
$$
By the construction, the measure $\mu$ is a symmetric measure for
the semigroup $\{T_t\}$ generated by the  process $X$:
$$
\E fT_tg=\E gT_tf,\quad t\geq 0;
$$
in particular,
$$
\E T_tf=\E f, \quad t\geq 0,
$$
i.e. $\mu$ is an invariant measure for $X$. The class
$\mathcal{G}=\mathcal{S}(\Re)$ is an algebra, invariant w.r.t. superpositions with
$C^\infty$-functions and dense in every
$L_p(\mu), p\geq 1$. In addition, thanks to the smoothness conditions and growth bounds imposed on coefficients $a, K$, the class $\mathcal{G}$
is invariant w.r.t. the semigroup $T_t$ and the generator $L$.  Define for $f, g\in \mathcal{G}$
$$
\Gamma(f,g)={1\over 2}\Big(L(fg)-fLg-gLf\Big), \quad
\Gamma_2(f,g)={1\over
2}\Big(L\Gamma(f,g)-\Gamma(Lf,g)-\Gamma(f,Lg)\Big).
$$
We will prove that \be\label{BE} \Gamma_2(f,f)\geq \alpha
\Gamma(f,f), \quad f\in \mathcal{G}, \ee then the required statement
would follow from the Bakry-Emery criterion \cite{BaEm}.

Straightforward calculations give
$$
\Gamma(f,g)=bf'g',
$$
$$\ba
2\Gamma_2(f,f)&=\Big(ab'+bb''-2a'b\Big)(f')^2-2bb'f'f''+2b^2(f'')^2\\&
=\left(ab'+bb''-2a'b-{(b')^2\over 2}\right)(f')^2  +\left({b'f'\over
\sqrt{2}}- bf''\sqrt{2}\right)^2\\&\geq
\left(ab'+bb''-2a'b-{(b')^2\over 2}\right)(f')^2. \ea
$$
Hence to prove (\ref{BE}) it is enough to show that \be\label{51}
2ab'+2bb''-4a'b-(b')^2\geq 4\alpha b. \ee Recall that
$$
v=\delta_\mu(K)=-K v_\mu-K',
$$
hence we can express the coefficients $a,b$ through the functions
$K$ and $v$:
$$
a=KK'-Kv, \quad b=K^2.
$$
Substituting these expressions into (\ref{51}), after some
transformations, which are straightforward but cumbersome and
therefore omitted, we re-write (\ref{51}) to the following form:
$$
K^3v'\geq \alpha K^2.
$$
The last inequality clearly holds true under (\ref{BE_K}). Hence,
applying the Bakry-Emery criterion, we get (\ref{LSI_K}) for every
$f\in \mathcal{S}(\Re)$.

If $K$ is bounded, then for every $f\in \mathcal{S}(\Re)$
(\ref{LSI}) holds true as a corollary of (\ref{LSI_K}).  It is a
standard procedure to approximate a given absolutely continuous $f$
such that $f, f'\in L_2(\mu)$ by a sequence of smooth compactly
supported $f_n$ in such a way that $f_n\to f$ and $f'_n\to f'$ in
$L_2(\mu)$; see e.g. the proof of Corollary 2.6.10 in
\cite{Bogachev}. Passing to the limit in (\ref{LSI}) for $f_n, n\geq
1$, we complete the proof. \qed

There is a wide choice for the pair of functions $v, K$ related by
(\ref{div}). Below we  give two versions of Theorem \ref{t3} which
correspond to particular choices of this pair. The first one arise
when one just takes  $v(x)=x-\langle\mu\rangle,$
$$
\langle\mu\rangle=\int_{\Re}y\, \mu(dy).
$$

\begin{prop}\label{p2} Let measure $\mu$ on $\Re$  have the first absolute moment  and have a positive continuous distribution density $p_\mu$. Denote
$$
\bar K_\mu(x)={1\over p_\mu(x)}\int_{x}^\infty
\Big(y-\langle\mu\rangle\Big)p_\mu(y)\, dy, \quad x\in \Re.
$$
The following statements hold true.

\begin{itemize}
  \item[I.] If $\inf_x\bar K_\mu(x)=\alpha>0$, then  for every smooth $f$ with a compact support
$$
\Ent f^{2} \leq \frac{2}{\alpha}\E\Big(\bar K_\mu f^{'}\Big)^2.$$
  \item[II.] If, in addition,  $\sup_x \bar K_\mu(x)=\beta<\infty$, then  for every absolutely continuous $f$ such that both $f$ and $f'$ are square integrable w.r.t. $\mu$,
$$
\Ent f^{2} \leq 2\bar c_\mu \E(f^{'})^2$$ with
$$
\bar c_\mu=\frac{\beta^2}{\alpha}.
$$
\end{itemize}

\end{prop}

In the second version of Theorem \ref{t3}, we choose $K$ in a  more
intrinsic way, namely, we take $K$ such that $\delta_\mu(K)=v$ with
\be\label{v_hat}v=\Phi^{-1}(F_\mu),\quad  F_\mu(x)=\mu((-\infty,
x]), \ee then $\mu\circ v^{-1}=\gamma$,
$\gamma\sim\mathcal{N}(0,1)$. Such a choice of the weight $v$ is
motivated by our intent to have
$$
\hat Z^v(\mu)=\hat Z(\gamma);
$$
that is, to make the order condition (\ref{LZ_v}) with $c=1$ as
precise as it is possible, i.e. to replace an inequality by an
identity. Because $\hat Z^v(\mu)=Z(\mu\circ v^{-1})$ identifies the
law of $\mu\circ v^{-1}$ uniquely, such an intent naturally leads to
the formula (\ref{v_hat}).

\begin{prop}\label{p3} Let measure $\mu$ on $\Re$ have a positive continuous distribution density $p_\mu$. Denote
$$
\hat K_\mu(x)={I(F_\mu(x))\over p_\mu(x)}.
$$
The following statements hold true.
\begin{itemize}
  \item[I.] For every smooth $f$ with a compact support,
\be\label{61} \Ent f^{2} \leq 2\E\Big(\hat K_\mu f^{'}\Big)^2.\ee
\item[II.] If, in addition,  $\hat K_\mu$ is bounded, then  for every absolutely continuous $f$ such that both $f$ and $f'$ are square integrable w.r.t. $\mu$,
$$
\Ent f^{2} \leq 2\hat c_\mu\E(f^{'})^2$$ with
$$
\hat c_\mu=\sup_x\Big(\hat K_\mu(x)\Big)^2.
$$
\end{itemize}
\end{prop}

\begin{rem} Define  the \emph{isoperimetric function} of the measure $\mu$ by
$$
I_\mu(p)=p_\mu(F_\mu^{-1}(p)), \quad p\in (0,1), \quad
I_\mu(0)=I_\mu(1)=1.
$$
 Then, clearly, the function $I$  defined by (\ref{iso}) equals $I_\gamma,\gamma \sim\mathcal{N}(0,1)$. The function $\hat K_\mu(x)$ above can be expressed as the ratio
$$
{I_\gamma(p)\over I_\mu(p)}\Big|_{p=F_\mu(x)},
$$
and under the conditions of Proposition \ref{p3} the function
$F_\mu$ gives a one-to-one correspondence between $(-\infty, \infty)$
and $(0,1)$. Hence the constant $\hat c_\mu$ above can be
alternatively expressed as
$$
\hat c_\mu=\left(\sup_{p\in (0,1)}{I_\gamma(p)\over
I_\mu(p)}\right)^2.
$$
\end{rem}

{\it Proofs of Proposition \ref{p2} and Proposition \ref{p3}.}
If  $v(x)=x-\langle\mu\rangle,$ we have $ \bar K_\mu v'=\bar  K_\mu,
$ and therefore the assumption $\inf \bar  K_\mu=\alpha>0$ made in
Proposition \ref{p2} implies the principal condition (\ref{BE_K}).
For the function $v$ defined by (\ref{v_hat}) and the function $\hat
K_\mu$, this condition takes even a more simple form because
straightforward calculation shows that
$$
\hat K_\mu v'=1.
$$
Hence one can expect that statements of  Proposition \ref{p2} and
Proposition \ref{p3} would follow from the version of the
Bakry-Emery criterion given in Theorem \ref{t3}. However, we can not
apply this theorem here directly, because of extra smoothness and
growth conditions on functions $K,a$, imposed therein. The strategy
of the proof will be the following: first, we consider a  family of
measures, which approximate $\mu$ properly and satisfy both
(\ref{BE_K}) for the respective pair of $K,v$, and extra smoothness
and growth conditions on respective functions $K,a$. Then, by
passing to a limit, we get respective  weighted log-Sobolev
inequality, i.e. prove statements I in Propositions \ref{p2},
\ref{p3}. Finally, using the same approximation procedure as in
the proof of Theorem \ref{t3} above, we extend the class of $f$ in the case
where the weight $K$ is bounded.

To shorten the exposition, we explain in details the way this
strategy is implemented for the proof of Proposition \ref{p3}, only.
The detailed proof of Proposition \ref{p2} is similar and omitted.
We also does not repeat the approximation arguments from the proof of Theorem \ref{t3}  above, and concentrate on the proof of (\ref{61}) for
smooth  compactly supported $f$.

Consider first the following auxiliary case: $p_\mu\in C^\infty$,
and for some $R>0$ \be\label{62} p_\mu(x)=\phi(x), \quad |x|\geq R.
\ee Then $v_\mu$  (which, let us recall, equals $p_\mu'/p_\mu$) and
$\hat K_\mu$ belong to $C^\infty$ and
$$
v_\mu(x)=-x, \quad \hat K_\mu(x)=1, \quad |x|\geq R.
$$
Then the functions $K=\hat K_\mu$ and $a$ defined by (\ref{ab})
satisfy the assumptions of Theorem \ref{t3}. Hence, applying Theorem
\ref{t3}, we get (\ref{61}).

Next, consider the general case. Fix some function $\chi\in
C^\infty$ taking values in $[0,1]$, such that $\chi(0)=0, \chi(x)=1,
x\geq 1$, and define
$$
\phi_{r, \delta}(x)=\phi(x)\Big(\delta+(1-\delta)\chi(|x|+r)\Big),
\quad x\in \Re;
$$
then every $\phi_{r, \delta}, r>0, \delta\geq 0$ belongs to
$C^\infty$. Denote
$$
M(r)=\int_{\Re}\phi_{r,0}(x)\, dx,
$$
then $M$  is a strictly decreasing function on $[0,\infty)$ and
$M(0)<1$. For a given $Q>0$, consider the restriction $p^Q_\mu$ of
$p_\mu$ to the segment $[-Q,Q]$, and assume that  $Q$ is large
enough for
$$
\int_{|x|>Q}p_\mu(x)\, dx<M(0).
$$
Then for every $\delta$ small enough there exists unique $r=r(Q,
\delta)>0$ such that
$$
\int_{\Re}\Big(p_\mu^Q(x)+\phi_{r, \delta}(x)\Big)\, dx=1.
$$
Take some non-negative  $\psi\in C^\infty$, supported in $[-1,1]$
and such that $\int_{\Re}\psi(x)\, dx=1$, and  consider the
probability measure $\mu_{Q,\eps}$ with the density
$$
p_{\mu_{Q,\delta}}(x)={1\over \delta}\int_{[-\delta,
\delta]}p_\mu^Q(y)\psi\left(x-y\over \delta\right)\, dy+\phi_{r,
\delta}(x).
$$
By the construction, every $\mu_{Q,\delta}$ has positive $C^\infty$
density and satisfy (\ref{62}) for some large $R$. Therefore,
(\ref{61}) holds true with $\mu_{Q,\delta}$ instead of $\mu$. It can
be seen easily that
$$
p_{\mu_{Q,\delta}}\to p_\mu, \quad K_{\mu_{Q,\delta}}\to K_\mu,
\quad \delta\to 0,  \quad Q\to \infty,
$$
uniformly on every finite segment. Passing to the limit, we obtain
(\ref{61}) for the initial measure $\mu$ and arbitrary  smooth and
compactly supported $f$.\qed

\section{Examples}\label{s4}

\begin{ex} Let $\mu$ on $\Re$  have a positive $C^1$-density $p_\mu$, such that for some $a, R>0$
\be\label{anti_lyap} v_\mu(x)x\geq -a x^2, \quad |x|>R \ee Let us
show that then condition $\inf \bar K_\mu>0$ from Proposition
\ref{p2} holds true. Changing the variables $x\mapsto x-\langle
\mu\rangle$, we can restrict ourselves to the case of $\langle
\mu\rangle =0$. Then we have for $x>R$
$$\ba
\bar K_\mu(x)&=\int_x^\infty y\exp({\log p_\mu(y)-\log p_\mu(x)})\,
dy\\&= \int_x^\infty y\exp\left(\int_x^y v_\mu(z)\, dz\right)\,
dy\geq \int_x^\infty y\exp\left(-a\int_x^y z\, dz\right)\, dy
\\&=e^{ax^2/2} \int_x^\infty ye^{-ay^2/2}\, dy=1/a.
\ea
$$
Similar relation holds true for $x<-R$; to see this, one should note
that
$$
\bar K_\mu=-{1\over p_\mu(x)}\int^{x}_{-\infty} yp_\mu(y)\, dy
$$
because $\mu$ is centered. Finally, because $p_\mu\in C^1$ is
positive, $\bar K_\mu$ has positive infimum over $[-R,R]$, which
completes the proof.

Similarly, if in addition for some $b>0$ \be\label{lyap}
v_\mu(x)x\leq -b x^2, \quad |x|>R, \ee then $\sup \bar
K_\mu<\infty$. Hence, by statement II of Proposition \ref{p2}, for a
measure $\mu$ satisfying (\ref{anti_lyap}) and (\ref{lyap}) the
log-Sobolev inequality holds true.

Note that (\ref{lyap}) is just the well known drift condition,
sufficient for the Poincar\'e inequality, e.g. Theorem 3.1 and
Remark 3.2 in \cite{CaGu}. However, various sufficient conditions
for the log-Sobolev inequality, available in the literature,
typically require additional assumptions on the \emph{curvature},
which in the current context equals $-v_\mu'$. Namely, the famous
Bakry-Emery condition (\cite{BaEm})) requires $-v_\mu'\geq
\delta>0$;  conditions by Wang (\cite{Wa}) and Cattiaux-Guillin
(\cite{CaGu}, Theorem 5.1) are more flexible, but still contain a
requirement that  the curvature is bounded from below, i.e. in our
case \be\label{curv_bound} -v_\mu'\geq \delta \ee with some
$\delta\in \Re$. The above condition (\ref{anti_lyap}) can be
understood as an ``integral'' version of (\ref{curv_bound}),
and it is easy to give an example of measure $\mu$ satisfying
(\ref{anti_lyap}) and (\ref{lyap}) such that (\ref{curv_bound})
fails.
\end{ex}

\begin{ex} Let $\gamma^3$ be a standard Gaussian measure on $\Re^3$, and $B_R$ be a ball of radius $R$, touching the origing and with the center located at the first basis vector $e_1$; that is, $B_R=B(Re_1,R)$.  Denote by $\gamma^{3,R}$ the measure $\gamma^3$ conditioned outside the ball $B_R$:
$$
\gamma^{3,R}(A)={\gamma^{3}(A\setminus B_R)\over
\gamma^{3}(\Re^3\setminus B_R)}.
$$
Consider a measure $\mu_R$ on $\Re$ which is a projection of
$\gamma^{3,R}$ on the first coordinate. We will show that there
exists some constant $\hat c$ such that uniformly by $R\geq 0$ the
constants $\hat c_{\mu}$ for the measures $\mu=\mu_R$ from
Proposition \ref{p3} are dominated by $\hat c$. This would yield
that for the family $\mu_R, R\geq 0$ the log-Sobolev inequality
holds true with uniformly bounded constants.

For a given $x\in [0, 2R]$,  the section of the ball $B_R$ by the
hyperplane $\{y=(y_1,y_2, y_3):y_1=x\}$, projected on the last two
coordinates, is the ball in $\Re^2$, centered at the origin and
having the radius
$$
r_R(x)=\sqrt{2Rx-x^2}.
$$
Define $r(x)=0$ for $x\not \in [0, 2R]$. Then we have for
$\mu=\mu_R$
$$
p_\mu(x)=C_R\phi(x)\psi_2\Big(r_R(x)\Big),
$$
where
$$
C_R=\Big(\gamma^{3}(\Re^3\setminus B_R)\Big)^{-1},
$$
$$
\psi_2(r)=\int_{\|y\|\geq r}{1\over
2\pi}e^{-(y_1^2+y_2^2)/2}\,dy_1dy_2={1\over
2\pi}\int_0^{2\pi}\int_r^\infty e^{-\rho^2/2}\rho\, d\rho
d\theta=e^{-r^2/2}.
$$
Consequently, \be\label{72} p_\mu(x)={C_R\over
\sqrt{2\pi}}\begin{cases}e^{-Rx},&x\in
[0,2R],\\e^{-x^2/2},&\hbox{otherwise}.\end{cases} \ee

To bound $\hat K_\mu(x)$ consider separately three cases.

\textbf{I.} $x<0$. Recall that $I'(p)=-\Phi^{-1}(p)$. Then for any
$c>1$ we have
$$
[I(c\Phi(x))]'=-\Phi^{-1}(c\Phi(x))c \phi(x)\leq
(-x)c\phi(x)=c\phi'(x)
$$
because  $\Phi^{-1}$ is an increasing  function. Clearly, both
$I(c\Phi(x))$ and $\phi(x)$ vanish as $x\to -\infty$, hence
\be\label{71} I(c\Phi(x))=\int_{-\infty}^x [I(c\Phi(y))]'\, dy\leq c
\int_{-\infty}^x \phi'(y)\, dy=c\phi(x), \quad x\leq \Phi^{-1}(1/c).
\ee Note that for $x<0$
$$
F_\mu(x)=C_R\Phi(x), \quad p_\mu(x)=C_R\phi(x),
$$
and $C_R>1$. In addition, the half-space $\{y=(y_1,y_2, y_3):y_1\leq
x\}$ is contained in $\Re^3\setminus B_R$, hence
$$
\Phi(x)=\gamma^3(\{y=(y_1,y_2, y_3):y_1\leq x\})\leq {1\over
C_R}\Leftrightarrow x\leq \Phi^{-1}\left({1\over C_R}\right),
$$
and we can apply (\ref{71}) to get
$$
\hat K_\mu(x)={I(C_R\Phi(x))\over C_R\phi(x)}\leq 1, \quad x<0.
$$

\textbf{II.} $x>2R$. In this case $1-F_\mu(x)=C_R(1-\Phi(x))$.
Recall that $I(p)=I(1-p)$ and $\Phi^{-1}\Big(1-\Phi(x)\Big)=-x$,
hence  we can use the same argument as  in the case \textbf{I} to
show that $\hat K_\mu(x)\leq 1$,  because for any $c>1$
$$
I(c(1-\Phi(x))=-c\int^{\infty}_x \Phi^{-1}\Big(c(1-\Phi(y)\Big)\,
dy\leq c \int^{\infty}_x y\phi'(y)\, dy=c\phi(x).
$$

\textbf{III.} $x\in [0, 2R]$. Recall that  there exists a constant
$c_*$ such that
$$
I(p)\leq c_*p\sqrt{\log\frac{1}{p}}, \quad p\in \left(0, {1\over
2}\right).
$$
One has
$$
C_R\gamma^3(\{y=(y_1,y_2, y_3):y_1>R\}))\leq 1-F_\mu(x)\leq
C_R\gamma^3(\{y=(y_1,y_2, y_3):y_1>0\}))<{1\over 2},
$$
hence we can write, using the identity $I(p)=I(1-p)$,
$$
\hat K_\mu(x)={I(1-F_\mu(x))\over p_\mu(x)}\leq c_*{1-F_\mu(x)\over
p_\mu(x)}\sqrt{\log\frac{1}{1-F_\mu(x)}}.
$$
Because $C_R>1$, we have
$$
\log\frac{1}{1-F_\mu(x)}\leq
\log\frac{1}{1-F_\mu(2R)}=\log\frac{1}{C_R(1-\Phi(2R))}\leq
\log\frac{1}{1-\Phi(2R)}\leq c^*(1+R)^2
$$
with some $c^*>2$. By (\ref{72}), we have
$$
{1-F_\mu(x)\over p_\mu(x)}=e^{Rx}\left(\int_x^{2R}e^{-Ry}\,
dy+\int_{2R}e^{-y^2/2}\, dy\right),
$$
and the right hand side term can be estimated either by
$$ e^{Rx}\int_x^{\infty}e^{-Ry}\, dy={1\over R}, $$
(when $R$ is large), or by
$$
e^{2R^2}\int_{0}e^{-y^2/2}\, dy={\sqrt{\pi\over 2}} e^{2R^2}
$$
(when $R$ is small). Then for any $R>0$ for $\mu=\mu_R$
$$
\hat c_\mu=\sup_x\hat K_\mu\leq \hat
c:=c_*c^*\sup_{Q>0}\min\left({1+Q\over Q}, {\sqrt{\pi\over
2}}(1+Q)e^{2Q^2}\right);
$$
for $R=0$ the measure $\mu$ just equals $\gamma$ and therefore $\hat
c_\mu=1$.

This example is motivated by the manuscript \cite{BCGM}, where the
problem of estimating of the Poincar\'e constant for a Gaussian
measure conditioned outside a ball is considered. One  approach
proposed therein is based on the decomposition of variance, and
requires an estimate for the Poincar\'e constant of one-dimensional
projection of the ``punctured'' Gaussian measure on the line which
contains the center of the ball. Such an estimate depend on the
position and the size  of the ball, see Lemma 4.7 in
\cite{BCGM}, and the case or a large ball touching the origin
relates the case (4) of that lemma. Our estimate for the log-Sobolev
constant implies that the Poincar\'e constant for $\mu$ is uniformly
bounded by $\hat c$, which drastically improves the bound $c
e^{R^2}$ from Lemma 4.7 \cite{BCGM}, statement (4). Heuristically,
the reason for this is the following. The measure $\mu$ contain
``cavities'',  which appear due to the ``puncturing'' procedure, and
if the ball is ``large'' and is located not so ``far from the origin'', then these
``cavities'' make the bounds for the Poincar\'e inequality obtained
via classic sufficient conditions to be very inaccurate. On the
other hand, the form of the weight  $\hat K_\mu$ in Proposition
\ref{p3} is highly adjusted to these ``cavities'', which makes
respective bounds more precise. We believe that similar calculations
can be made in a general setting, i.e. for arbitrary $d\geq 2$ and
arbitrary  position and size of the ball; this is a subject of a
further research.
\end{ex}

\subsection*{Acknowledgements} The paper was  prepared partially during a visit of authors to the
University of Potsdam; the authors express their deep gratitude to
the University of Potsdam and personally to Sylvie Roelly for the
hospitality. The authors are  also grateful to Emmanuel Boissard for
fruitful discussions concerning the  manuscript \cite{BCGM}.

\end{document}